\documentclass[a4paper,11pt]{article}
\usepackage[utf8]{inputenc}
\usepackage[T1]{fontenc}
\usepackage[english]{babel}
\usepackage[margin=2.8cm]{geometry}
\usepackage{mathpazo}
\usepackage[scaled=.95]{helvet}
\usepackage{courier}
\usepackage{float}
\usepackage{amsmath,amsfonts,amssymb,amsthm}
\usepackage{graphicx}
\usepackage[square,sort,comma,numbers]{natbib}
\usepackage{color}
\usepackage[onehalfspacing]{setspace}

\newtheorem{theo}{Theorem}[section]

\usepackage[draft]{fixme}
\fxusetheme{color}
\FXRegisterAuthor{h}{ah}{H}
\FXRegisterAuthor{p}{ap}{P}
\FXRegisterAuthor{a}{aa}{A}

\usepackage{hyperref}
\hypersetup{
pdfpagemode=UseOutlines,      
pdfstartview=Fit,             
pdffitwindow=true,            
pdfpagelayout=TwoColumnsRight,
pdftoolbar=true,              
pdfmenubar=true,              
bookmarksopen=false,          
bookmarksnumbered=true,       
colorlinks=true,              
pdfauthor={Ton nom},     
pdftitle={Titre PDF},    
pdfcreator=PDFLaTeX,          %
pdfproducer=PDFLaTeX,         %
linkcolor=blue,               
urlcolor=blue,                
anchorcolor=blue,            
citecolor=blue,              
frenchlinks=true,             
pdfborder={0 0 0}             
}

\begin{document}
\title{On the rates of convergence of Parallelized Averaged Stochastic Gradient Algorithms}
\author{Antoine Godichon-Baggioni$^*$, Sofiane Saadane$^{**}$\\ 
 $^*$ Laboratoire de Mathématiques de l'INSA de Rouen, \\
INSA de Rouen, 76800 Saint-Etienne-du-Rouvray, France \\
$^{**}$ Institut de Mathématiques de Toulouse,\\
INSA de Toulouse, 31000 Toulouse, France \\
 email: $^*$ antoine.godichon@insa-rouen.fr, $^{**}$ saadane@insa-toulouse.fr, 
} 
\date{}
\maketitle

\begin{abstract}
The growing interest for high dimensional and functional data analysis led in the last decade to an important research developing a consequent amount of techniques. Parallelized algorithms, which consist in distributing and treat the data into different machines, for example, are a good answer to deal with large samples taking values in high dimensional spaces. We introduce here a parallelized averaged stochastic gradient algorithm, which enables to treat efficiently and recursively the data, and so, without taking care if the distribution of the data into the machines is uniform. The rate of convergence in quadratic mean as well as the asymptotic normality of the parallelized estimates are given, for strongly and locally strongly convex objectives.
\end{abstract}

\medskip


\section{Introduction}

The growing interest for high dimensional and functional data analysis led in the last decade to many research papers developing a consequent amount of techniques. Data that have to face statisticians can now be extremely large so that it creates a need for economic calculation techniques. Parallelized algorithms are now a good answer to this challenge and authors are using various techniques and software (multicore processors for example). For instance, \citep{zinkevich2010parallelized} deals with gradient descent for least square type functions (strongly convex objective functions) by using a global averaging technique meaning that each machine carries out a gradient descent and the final outcome of the procedure is an averaging of all these results. In a recent work, \citep{recht2013parallel} also deals with stochastic gradient by proposing a rewriting procedure where each processor can rewrite the data of an other. This paper offers good numerical results proving its efficiency and showing in particular that the rewriting technique is sparse (meaning it does affect a too much important part of the data). We can also point out the work of \citep{liu2015asynchronous} where two implementations of stochastic gradient are coupled : one is over a computer network and the other one is on a shared memory system. Finally, \citep{bianchi2013performance} introduces a parallelized averaged stochastic gradient algorithm and establishes some convergence properties, such as the asymptotic normality. Those four examples show in particular that the literature assumes (most of the time) that each machine, when working with a parallelized algorithm, receive the same amount of data. The reality of course is far from this setup. Indeed, in many practical cases, data are acquired and treated by different machines (see the nice example of software architecture given in \citep{biau2014online}). For this reason, it is interesting to investigate a parallelized algorithm when the distribution of the data into the machines is not uniform. Moreover, data are often acquired sequentially, then, it is important to have an algorithm which enables to simply update the estimates. We so focus on the parallelization of averaged stochastic gradient algorithms.

\medskip
\noindent Stochastic Gradient Descents (SGD for short) are usually used for estimating the minimizer of a convex function and are commonly fast, do not need to store all the data into memory, and are recursive, which enables to simply update the estimates when the data arrive sequentially (\citep{robbins1951,Duf97,KY03}). In order to improve the convergence, \citep{ruppert1988efficient} and \citep{PolyakJud92} introduced the Averaged Stochastic Gradient Descent (ASGD for short). Among the studied cases, two of them attract a lot of attention : globally strongly convex objectives (\citep{bach2013non,ghadimi2012optimal}), and locally strongly convex ones (\citep{pelletier1998almost,Pel00,GB2016,GB2017}). Indeed, in those cases the theoretical study can be pushed very far due to the quite nice structure of the objective function. 

\medskip
\noindent In this paper, we introduce Parallelized Stochastic Gradient (PASG) algorithm, which consists in running $p$ samples of ASGD with sample sizes $n_i$, $i\in\lbrace 1,...,p\rbrace$. After a run, the results are centralized using an averaging step, \textit{i.e.} taking the arithmetic mean of all the estimates obtained with each SGD (or equivalently taking the weighted mean of those obtained with each ASGD). The interest of this procedure is its ability to deal with large samples not necessarily with the same size (we can suppose $n_i \neq n_j$ for any $i\neq j$). We then establish the efficiency of the algorithms by proving that they have a quadratic convergence rate of $O \left(\frac{1}{n} \right)$, which is the optimal one for stochastic approximation. In a second time, we establish the asymptotic normality of the estimates and see that it has an optimal asymptotic variance \citep{Pel00}.

\medskip
\noindent The paper is organized as follows: Section \ref{secdef}, some recalls on SGD and its averaged version are done and the general framework as well as the PASG algorithm are introduced. The two contexts (globally and locally convex objective) are introduced in Section \ref{secconv} as well as the rate of convergence in quadratic mean and the asymptotic normality of the estimates obtained with the PASG algortihm. Section \ref{secnum}, a simulation study for the estimation of the geometric median shows the efficiency of the method. Finally, the proofs are postponed in Section~\ref{secproof}.

\section{Framework and algorithms}\label{secdef}
\subsection{General framework and usual averaged stochastic gradient algorithm}
Let $X$ be a random variable taking values in a space $\mathcal{X}$, and let $H$ be a separable Hilbert space, not necessarily of finite dimension, such as $\mathbb{R}^{d}$ or $L^{2}(I)$, for some close interval $I \subset \mathbb{R}$. In what follows, we denote by $\langle .,. \rangle$ its inner product, and by $\| . \|$ the associated norm. Let $G: H \longrightarrow \mathbb{R}$ be the function we would like to minimize, defined for all $h \in H$ by
\begin{equation}
G(h) := \mathbb{E}\left[ g \left( X,h \right) \right],
\end{equation}
where $g: \mathcal{X}\times H \longrightarrow \mathbb{R}$. We consider from now that the functional $G$ is convex, and that for almost every $x \in \mathcal{X}$, the functional $g\left( x,. \right)$ is Fréchet-differentiable 
for the second variable and we denote by $\nabla_{h}g \left( x,. \right)$ its gradient. Let $X_{1},...,X_{n},...$ be random variables with the same law as $X$, the stochastic gradient descent (SGD for short) is defined recursively for all $n \geq 1$ by (\citep{robbins1951})
\begin{equation}
m_{n+1} = m_{n} - \gamma_{n} \nabla_{h} g \left( X_{n+1} , m_{n} \right) ,
\end{equation}
with $m_{1}$ bounded, and a step sequence $\left( \gamma_{n} \right)_{n \geq 1}$ of the form $ \gamma_{n} := c_{\gamma}n^{-\alpha}$, with $c_{\gamma}> 0$ and $\alpha \in \left( \frac{1}{2}, 1 \right)$. Remark that it is possible to take a step sequence of the form $\gamma_{n} = \frac{c}{n}$, but it necessitates to have some information on the smallest eigenvalue of the Hessian of the functional $G$ at $m$ (\citep{pelletier1998almost}). In order to improve the convergence, \citep{ruppert1988efficient} (see also \citep{PolyakJud92} for first results) introduced the averaged stochastic gradient descent (ASGD for short), defined recursively for all $n \geq 1$ by 
\begin{equation}
\overline{m}_{n+1} = \overline{m}_{n} + \frac{1}{n+1}\left( m_{n+1} - \overline{m}_{n} \right),
\end{equation}
with $\overline{m}_{1} = m_{1}$. We speak about averaging since it can be written as
\begin{equation}
\overline{m}_{n} = \frac{1}{n}\sum_{k=1}^{n} m_{k} .
\end{equation}

\subsection{The parallelized averaged stochastic gradient algorithm}
We consider from now a set $\left\lbrace 1,...,p \right\rbrace$ of machines. The data are spread over the machines, \textit{i.e.} each entity $i$ receives sequentially a sequence of independent random variables $X_{i,1},...,X_{i,k},...$. Then, each entity $i=1,...,p$ will compute the SGD and its averaged version defined recursively for all $k \geq 1$ by 
\begin{center}
$\left\{
\begin{array}{rl}
 m_{i,k+1} &= m_{i,k} - \gamma_{k}\nabla_{h} g \left( X_{i,k+1} , m_{i,k} \right) , \\
 \overline{m}_{i,k+1} &= \overline{m}_{i,k} + \frac{1}{k+1} \left( m_{i,k+1} - \overline{m}_{i,k} \right) ,  
\end{array}
\right.$
\end{center}
with $m_{1,1}=\overline{m}_{1,1},...,m_{p,1}= \overline{m}_{p,1}$ bounded and $\left( \gamma_{k} \right)_{k\geq 1}$ a step sequence of the form $\gamma_k = c_{\gamma}k^{-\alpha}$, with $c_{\gamma}>0$ and $\alpha \in \left( \frac{1}{2},1 \right)$. Let $n = n_{1}+...+n_{p} $, the parallelized averaged stochastic gradient estimate at time $n$ is defined by
\begin{equation}
\widehat{m}_{n} := \frac{1}{\sum_{i=1}^{p}n_{i}}\sum_{i=1}^{p}n_{i} \overline{m}_{i,n_{i}} ,
\end{equation} 
which can be written as
\begin{equation}
\widehat{m}_{n} = \frac{1}{\sum_{i=1}^{p}n_{i}}\sum_{i=1}^{p}\sum_{k=1}^{n_{i}}m_{i,k} .
\end{equation}
Moreover, this algorithm can be written recursively. Indeed, setting $n' = n_{1}'+...+n_{p}' $, such that for all $i=1,...,p$ we have $n_{i}' \geq n_{i}$, 
\begin{equation}
\widehat{m}_{n'} = \frac{\sum_{i=1}^{p}n_{i}}{\sum_{i=1}^{p}n_{i}'} \widehat{m}_{n} + \frac{1}{\sum_{i=1}^{p}n_{i}'}\sum_{i=1}^{p} \left( n_{i}'\overline{m}_{i,n_{i}'} - n_{i}\overline{m}_{i,n_{i}}\right) .
\end{equation}

\section{Convergence results}\label{secconv}
\subsection{Strongly convex objective}
We now introduce sufficient conditions which ensures the convergence of stochastic gradient algorithms and of the PASG-algorithm when the functional $G$ is strongly convex. 
\begin{itemize}
\item[\textbf{(H1)}] The functional $G$ is differentiable and denoting by $\Phi$ its gradient, there exists $m \in H$ such that
\[
\Phi (m) := \nabla G (m) = 0 .
\]
\item[\textbf{(H2)}] The functional $G$ is twice continuously differentiable almost everywhere and for all positive constant $A$, there is a positive constant $C_{A}$ such that for all $h \in \mathcal{B}\left( m , A \right)$,
\[
\left\| \Gamma_{h} \right\|_{op} \leq C_{A} ,
\]
where $\Gamma_{h}$ is the Hessian of the functional $G$ at $h$ and $ \left\| . \right\|_{op}$ is the usual spectral norm for linear operators.\\

\item[\textbf{(H3)}] There is a positive constant $C_{m}$ such that for all $h \in H$,
\[
\left\| \nabla G (h) - \Gamma_{m}(h-m) \right\| \leq C_{m} \left\| h-m \right\|^{2}.
\]
\item[\textbf{(H4)}]  There are positive constants $L_{1}, L_{2}$ such that for all $h \in H$,
\begin{align*}
& \mathbb{E}\left[ \left\| \nabla_{h}g \left( X,h \right) \right\|^{2} \right] \leq L_{1} \left( 1 + \left\| h-m \right\|^{2} \right) , \\
& \mathbb{E}\left[ \left\| \nabla_{h}g \left( X,h \right) \right\|^{4} \right] \leq L_{2} \left( 1 + \left\| h-m \right\|^{4} \right) . 
\end{align*}
\item[\textbf{(H5)}] The functional $G$ is $\mu$-strongly convex: for all $h,h' \in H$,
\[
G(h) \geq G(h') + \left\langle \nabla G (h) , h'-h \right\rangle + \mu \left\| h-h' \right\|^{2} .
\] 
\end{itemize}
We now make some comments on the assumptions: \textbf{(H1)} simply states the existence of a local minimum which is a necessary condition for our work. Assumptions \textbf{(H2)} to \textbf{(H5)} are smoothness properties stating that $G$ is $\mu$-strongly convex, coercive and have at most quadratic growth. \textbf{(H5)} is still standard and is the most favorable case when dealing with convex optimization problems, leading to the best possible achievable rates. Remark that the literature is very large on the rate of convergence of stochastic gradient algorithms in the case of strongly convex objective (see \citep{bach2013non} among others) and one can check that under assumptions \textbf{(H1)} to \textbf{(H5)}, there are positive constants $C_{1},C_{2}$ such that for all $n \geq 1$,
\begin{align}
& \mathbb{E}\left[ \left\| m_{n} - m \right\|^{2} \right] \leq \frac{C_{1}}{n^{\alpha}}, \label{vitl2}\\
& \mathbb{E}\left[ \left\| m_{n} - m \right\|^{4} \right] \leq \frac{C_{2}}{n^{2\alpha}} \label{vitl4}.
\end{align}
The following theorem gives the rates of convergence in quadratic mean of the PASG-algorithm.

\begin{theo}\label{theostr}
Suppose assumptions \textbf{(H1)} to \textbf{(H5)} hold. Then, for all $n=\sum_{i=1}^{p}n_{i}$,
\begin{align*}
\mathbb{E}\left[ \left\| \widehat{m}_{n} - m \right\|^{2} \right] & \leq \frac{L_{1}\lambda_{\min}^{-2}}{\sum_{i=1}^{p}n_{i}} + \sum_{j=1}^{5}\lambda_{\min}^{-2}A_{j,p,n}^2 + \sum_{j,j'=1, j \neq j'}^{6}\lambda_{\min}^{-2}A_{j,p,n}A_{j',p,n},
\end{align*}
where $\lambda_{\min} \geq \mu $ is the smallest (or limit inf for infinite dimensional spaces) eigenvalue of $\Gamma_{m}$ and
\begin{eqnarray*}
 A_{1,p,n}^2 = A_{2,p,n}^2 := \frac{p^{2}C_{1}c_{\gamma}^{-2}}{\left( \sum_{i=1}^{p}n_{i} \right)^{2}}, \quad \quad  &  A_{3,p,n}^2 := \frac{4p^{2-\alpha}\alpha c_{\gamma}^{-2}C_{1}}{\left( \sum_{i=1}^{p}n_{i} \right)^{2-\alpha}}, \quad \quad 
&  A_{4,p,n}^2 := \frac{C_{m}^{2}C_{2}\left( 1-\alpha \right)^{-2}p^{2\alpha}}{\left( \sum_{i=1}^{p}n_{i} \right)^{2\alpha}}, \\ 
 A_{5,p,n}^2 := \frac{L_{1}C_{1}\left( 1-\alpha \right)^{-1}p^{\alpha}}{\left( \sum_{i=1}^{p}n_{i} \right)^{1+\alpha}},   \quad \quad & A_{6,p,n}^2 := \frac{L_{1}}{\sum_{i=1}^{p}n_{i}}. \quad \quad &
\end{eqnarray*}
More precisely, we have
\[
\mathbb{E}\left[ \left\| \widehat{m}_{n} - m \right\|^{2} \right]  \leq \frac{L_{1}\lambda_{\min}^{-2}}{\sum_{i=1}^{p}n_{i}} + o \left( \frac{1}{\sum_{i=1}^{p}n_{i}} \right) .
\]
\end{theo}

\noindent\textbf{Remark:} 
\begin{itemize}
\item One can note that the rate of convergence is the optimal one for strongly convex function (see \citep{nemirovsky1983problem} for instance) and that the choice $\alpha\in(1/2,1)$ is crucial to obtain this bound. Indeed, when $\alpha\in(0,1/2)$, the result is different mainly because remainder terms play a preponderant role while our choice is justified by the central limit theorem. Indeed, the rate of convergence can be considered as optimal in our case since it perfectly reflects the asymptotic normality (see Theorem \ref{tlcstr}).
\item Investigating the case $\alpha =1$ is a tricky question discussed before, and not necessary since optimality is already obtained. However, we point out the work of (\citep{bach2013non} and \citep{GPS16}) where the specificity of this case is discussed with accurate computations. 
\item Remark that the remainders terms are negligible since $p = o \left( n^{ \max \left\lbrace \frac{2\alpha -1}{2\alpha} , \frac{1-\alpha}{2-\alpha}\right\rbrace} \right)$. For example, for $\alpha = \frac{2}{3}$, they are negligible since $p = o \left( \sqrt{n} \right)$.
\item Among the classical examples, one can think about least-square regression, where the objective function is of the form $\mathbb{E}\left[ \left( \left\langle X,h \right\rangle - Y \right)^{2} \right]$, for $X \in \mathbb{R}^{d}$ and $Y \in \mathbb{R}$ (see \citep{dieuleveut2016harder,cohen2017projected} for instance).
\end{itemize}

In order to establish a Central Limit Theorem, let us now introduce a new assumption:
\begin{itemize}
\item[\textbf{(H6)}] Let $\|.\|_{F}$ be the Frobenius norm for linear operators,
\[
\lim_{h \to m} \left\| \mathbb{E}\left[ \nabla_{h} g \left( X , m \right) \otimes \nabla_{h}g \left( X,m \right) \right] - \mathbb{E}\left[ \nabla_{h} g \left( X , h \right) \otimes \nabla_{h} g\left( X,h \right) \right] \right\|_{F} = 0 ,
\]
where for all $h,h',h'' \in H$, $h\otimes h'(h'')=\langle h,h''\rangle h'$.
\end{itemize}
We can now give the asymptotic normality of $\left( \widehat{m}_{n} \right)$.
\begin{theo}\label{tlcstr}
Suppose assumptions \textbf{(H1)} to \textbf{(H6)} hold. Then, let $n=\sum_{i=1}^{p}n_{i}$,
\[
\lim_{n\to\infty}\sqrt{\sum_{i=1}^{p}n_{i}}\left( \widehat{m}_{n}-m \right) \sim \mathcal{N}\left( 0 , \Gamma_{m}^{-1}\Sigma \Gamma_{m}^{-1} \right) ,
\]
where 
\[
\Sigma:=\mathbb{E}\left[ \nabla_{h}g \left( X,m \right) \otimes \nabla_{h}g\left( X,m\right) \right].
\]
\end{theo}

\subsection{Locally strongly convex objective}
We now focus on the framework introduced by \citep{GB2016} and \citep{GB2017} when $G$ is only locally strongly convex:
\begin{itemize}
\item[\textbf{(H7)}] There exists a positive constant $\epsilon$ such that for all $h \in \mathcal{B}\left( m , \epsilon \right)$, there is a basis of $H$ composed of eigenvectors of $\Gamma_{h}$. Moreover,  let us denote by $\lambda_{\min}$ the limit inf of the eigenvalues of $\Gamma_{m}$, then $\lambda_{\min}$ is positive. Finally, for all $h \in \mathcal{B}\left( m , \epsilon \right)$, and for all eigenvalue $\lambda_{h}$ of $\Gamma_{h}$, we have $\lambda_{h} \geq \frac{\lambda_{\min}}{2} > 0$.
\item[\textbf{(H8)}]
For all integer $q$, there is a positive constant $L_{q}$ such that for all $h \in H$,
\[
\mathbb{E}\left[ \left\| \nabla_{h}g \left( X,h \right) \right\|^{2q} \right] \leq L_{q} \left( 1 + \left\| h-m \right\|^{2q} \right) . 
\]
\end{itemize}
The main difference with previous framework is that we just have to assume the local strong convexity of the functional we would like to minimize, and in return, we have to assume the existence of the $q$-th moments of the gradient. Note that under assumptions \textbf{(H1)} to \textbf{(H4)} and \textbf{(H7)}, it was proven that for all positive constant $\delta$,
\begin{equation}\label{vitas}
\left\| m_{n} - m \right\|^{2} = o \left( \frac{(\ln n)^{\delta}}{n^{\alpha}}\right).
\end{equation}
Moreover, suppose assumption \textbf{(H8)} holds too, it was proven that for all positive integer $p$, there is a positive constant $C_{p}$ such that for all $n \geq 1$,
\begin{equation}\label{vitlp}
\mathbb{E}\left[ \left\| m_{n} - m \right\|^{2p}\right] \leq \frac{C_{p}}{n^{p\alpha}}.
\end{equation}
Then, we can now give the rate of convergence in quadratic mean of the PASG-algorithm for locally strongly convex objectives.
\begin{theo}\label{theoloc}
Suppose assumptions \textbf{(H1)} to \textbf{(H3)} and \textbf{(H7)}, \textbf{(H8)} hold. Then, for all $n=\sum_{i=1}^{p}n_{i}$,
\begin{align*}
\mathbb{E}\left[ \left\| \widehat{m}_{n} - m \right\|^{2} \right] & \leq \frac{L_{1}\lambda_{\min}^{-2}}{\sum_{i=1}^{p}n_{i}} + \sum_{j=1}^{5}\lambda_{\min}^{-2}A_{j,p,n}^2 + \sum_{j,j'=1, j \neq j'}^{6}\lambda_{\min}^{-2}A_{j,p,n}A_{j',p,n},
\end{align*}
where
\begin{eqnarray*}
 A_{1,p,n}^2 = A_{2,p,n}^2 := \frac{p^{2}C_{1}c_{\gamma}^{-2}}{\left( \sum_{i=1}^{p}n_{i} \right)^{2}}, \quad \quad  &  A_{3,p,n}^2 := \frac{4p^{2-\alpha}\alpha c_{\gamma}^{-2}C_{1}}{\left( \sum_{i=1}^{p}n_{i} \right)^{2-\alpha}}, \quad \quad 
&  A_{4,p,n}^2 := \frac{C_{m}^{2}C_{2}\left( 1-\alpha \right)^{-2}p^{2\alpha}}{\left( \sum_{i=1}^{p}n_{i} \right)^{2\alpha}}, \\ 
 A_{5,p,n}^2 := \frac{L_{1}C_{1}\left( 1-\alpha \right)^{-1}p^{\alpha}}{\left( \sum_{i=1}^{p}n_{i} \right)^{1+\alpha}},   \quad \quad & A_{6,p,n}^2 := \frac{L_{1}}{\sum_{i=1}^{p}n_{i}}. \quad \quad &
\end{eqnarray*}
More precisely, we have
\[
\mathbb{E}\left[ \left\| \widehat{m}_{n} - m \right\|^{2} \right]  \leq \frac{L_{1}\lambda_{\min}^{-2}}{\sum_{i=1}^{p}n_{i}} + o \left( \frac{1}{\sum_{i=1}^{p}n_{i}} \right) .
\]
\end{theo}

Finally, we also establish the asymptotic normality of $\left( \widehat{m}_{n} \right)$.
\begin{theo}\label{tlcloc}
Suppose assumptions \textbf{(H1)} to \textbf{(H4)} and \textbf{(H6)}, \textbf{(H7)} hold. Then, let $n=\sum_{i=1}^{p}n_{i}$,
\[
\lim_{n\to\infty}\sqrt{\sum_{i=1}^{p}n_{i}}\left( \widehat{m}_{n}-m \right) \sim \mathcal{N}\left( 0 , \Gamma_{m}^{-1}\Sigma \Gamma_{m}^{-1} \right) ,
\]
where 
\[
\Sigma:=\mathbb{E}\left[ \nabla_{h}g \left( X,m \right) \otimes \nabla_{h}g\left( X,m\right) \right].
\]
\end{theo}

\noindent\textbf{Remark:}
\begin{itemize}
\item Among the classical examples, one can think about logistic regression \citep{bach2014adaptivity}, which leads to minimize $\mathbb{E}\left[ \ln \left( 1 + \exp \left( -Y \left\langle X,h \right\rangle \right)\right) \right] $, where $Y \in \lbrace -1,1 \rbrace$ and $X \in \mathbb{R}^{d}$. One can also think about the estimation of the geometric median (see \citep{Hal48,HC,CCG2015} among others), where the objective function is $ \mathbb{E}\left[ \left\| X - h \right\| - \left\| X \right\| \right]$.
\end{itemize}

\section{Numerical experiments}\label{secnum}
\begin{figure}[H]
\begin{center}
\includegraphics[scale=0.45]{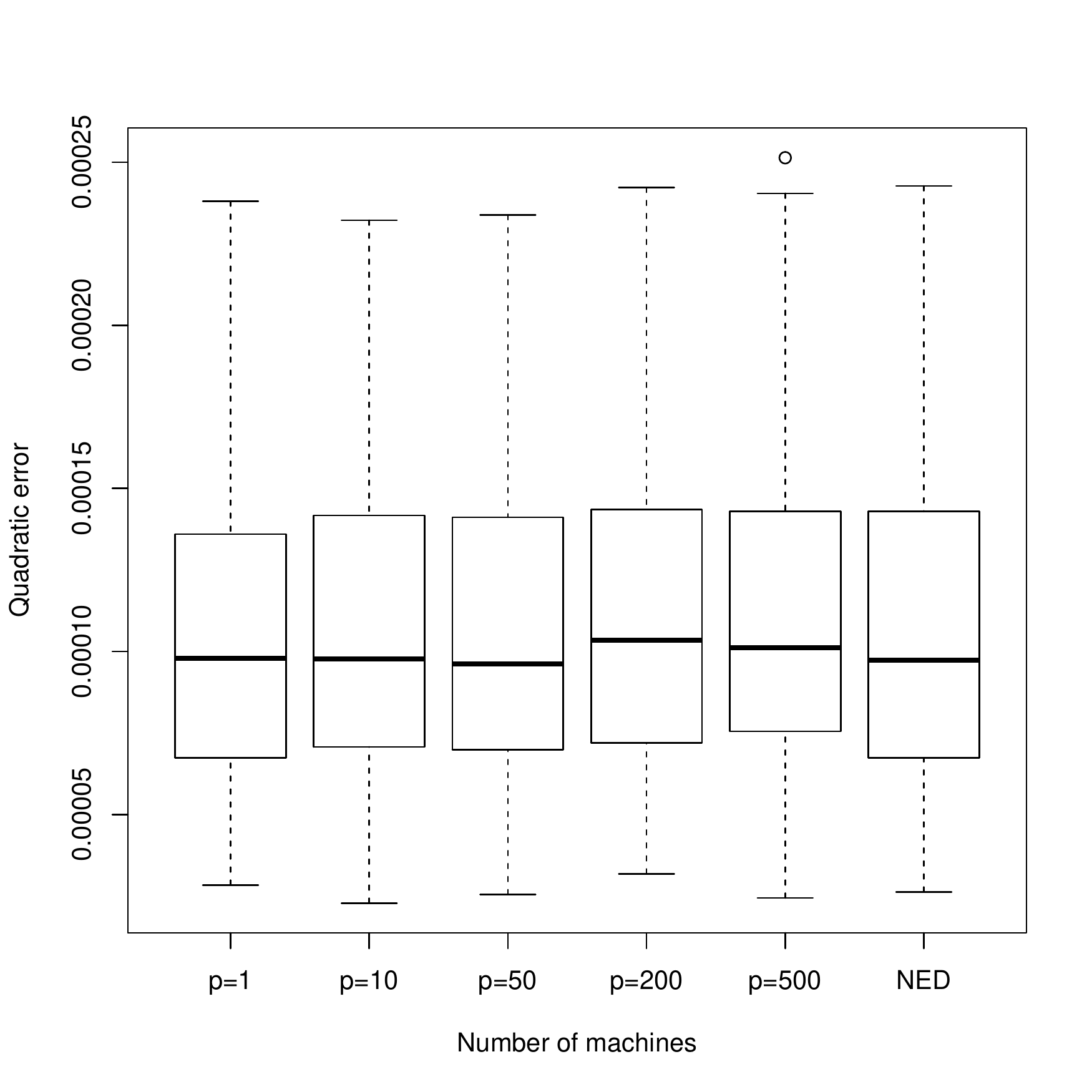}
\caption{Quadratic error obtained with the PASGD for a sample size $n=10^5$ and for different numbers of machines ($p=1,10,50,200,500$) and for a Non-Equal Distribution ("NED") for $p=10$. \label{n105}}
\end{center}
\end{figure}

In this section, in order to illustrate our theoretical work, we propose some numerical experiments. We consider from now a Gaussian vector $X$ taking values in $\mathbb{R}^{d}$, with mean $\mathbb{E}[X]=0$ and variance $\mathbb{E}[ X \otimes X ] = I_{d}$. We focus on the estimation of the Geometric Median of $X$, which is defined by \citep{Hal48,Kem87}
\[
m:= \arg \min_{h \in \mathbb{R}^{d}} \mathbb{E}\left[ \left\| X-h \right\| - \left\| X \right\| \right] .
\]
Note that in this case, $m=0$, and the SGD is defined for all $i=1,\ldots ,p$ and $k \geq 1$ by \citep{HC}
\begin{align*}
m_{i,k+1} = m_{i,k} + \gamma_{k} \frac{X_{i,k+1} - m_{i,k}}{\left\| X_{i,k} - m_{i,k}\right\|}.
\end{align*}
We now consider a step sequence $\gamma_{k}=k^{-2/3}$ and numbers of machines equal to $1,10,50,200,500$, and we assume that data are uniformly distributed into the different machines. Note that in this case, our algorithm is the same as the one introduced by \citep{bianchi2013performance}. Moreover, when $p=1$, this corresponds to the usual ASG algorithm. Finally, we also consider a case with $p=10$ but where we have a Non-Equal Distribution ("NED" for short) between the machines, and the following vector gives the percentage of data per machine:
\[
v_{p} = \left( 0.05, 0.45, 1.5,   3, 8, 10, 10, 17,20, 30 \right)
\]
In Figure \ref{n105}, one can see that the number of used machines (with $p  \ll \sqrt{n}$) does not seem to have a strong impact on the quadratic error of the estimates, which tends to confirm the results given by Theorems \ref{theostr} and \ref{theoloc}. Moreover, it seems to confirm that taking into account the numbers of data per machine during the parallelization step enables to take care of estimation coming from strongly different and inhomogeneous sources.  

\begin{figure}[H]
\begin{center}
\includegraphics[scale=0.45]{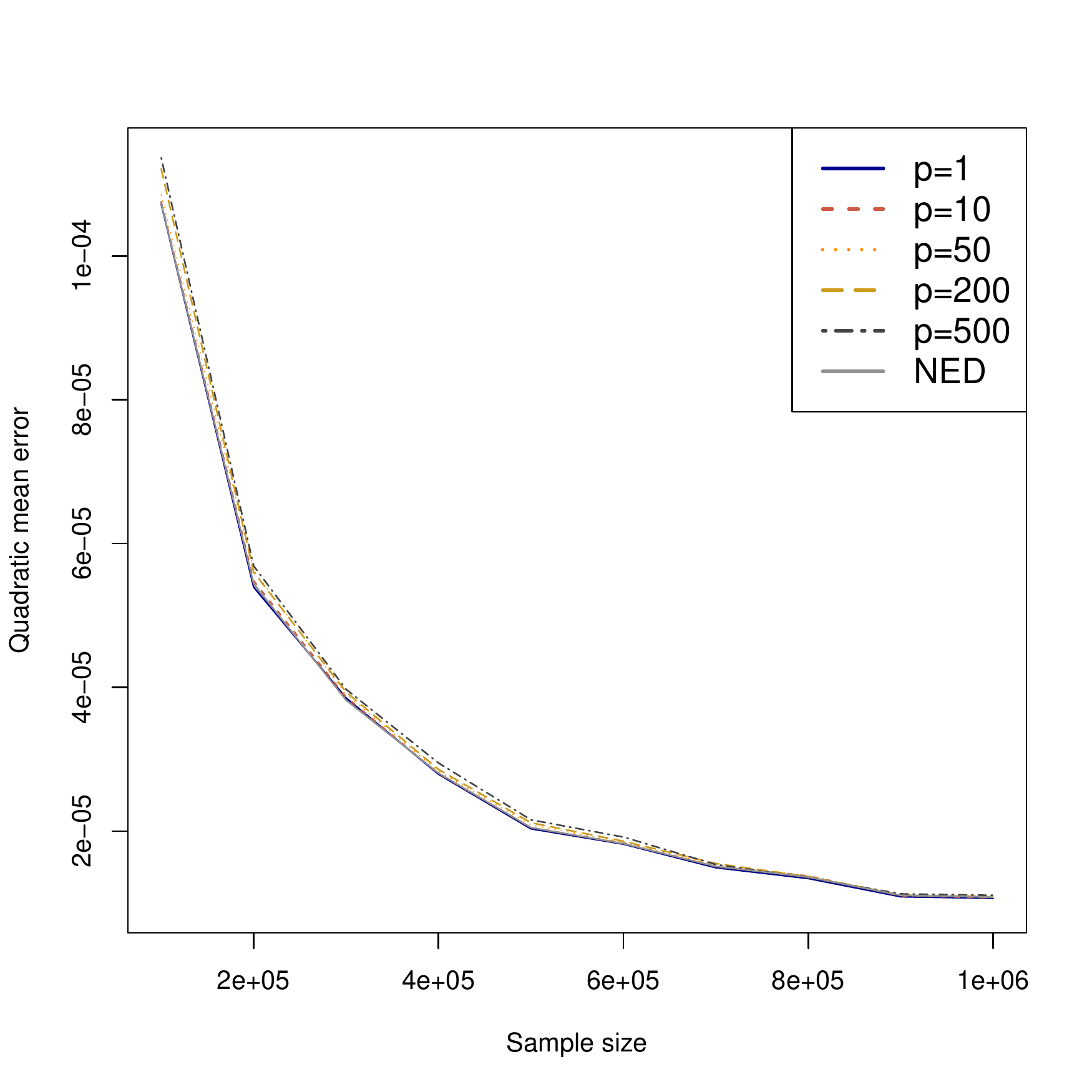}
\caption{Evolution of the quadratic mean error compare to the sample size $n$ for different numbers of machines ($p=1,10,50,200,500$) and for a Non-Equal Distribution ("NED") for $p=10$.\label{figev}}
\end{center}
\end{figure}

Figure \ref{figev} tends to confirm Theorems \ref{tlcstr} and \ref{tlcloc}, i.e it tends to confirm that the asymptotic behavior of the PASG algorithm does not depend on the number of machines or on the homogeneity of the data distribution.

\begin{figure}[H]
\begin{center}
\includegraphics[scale=0.45]{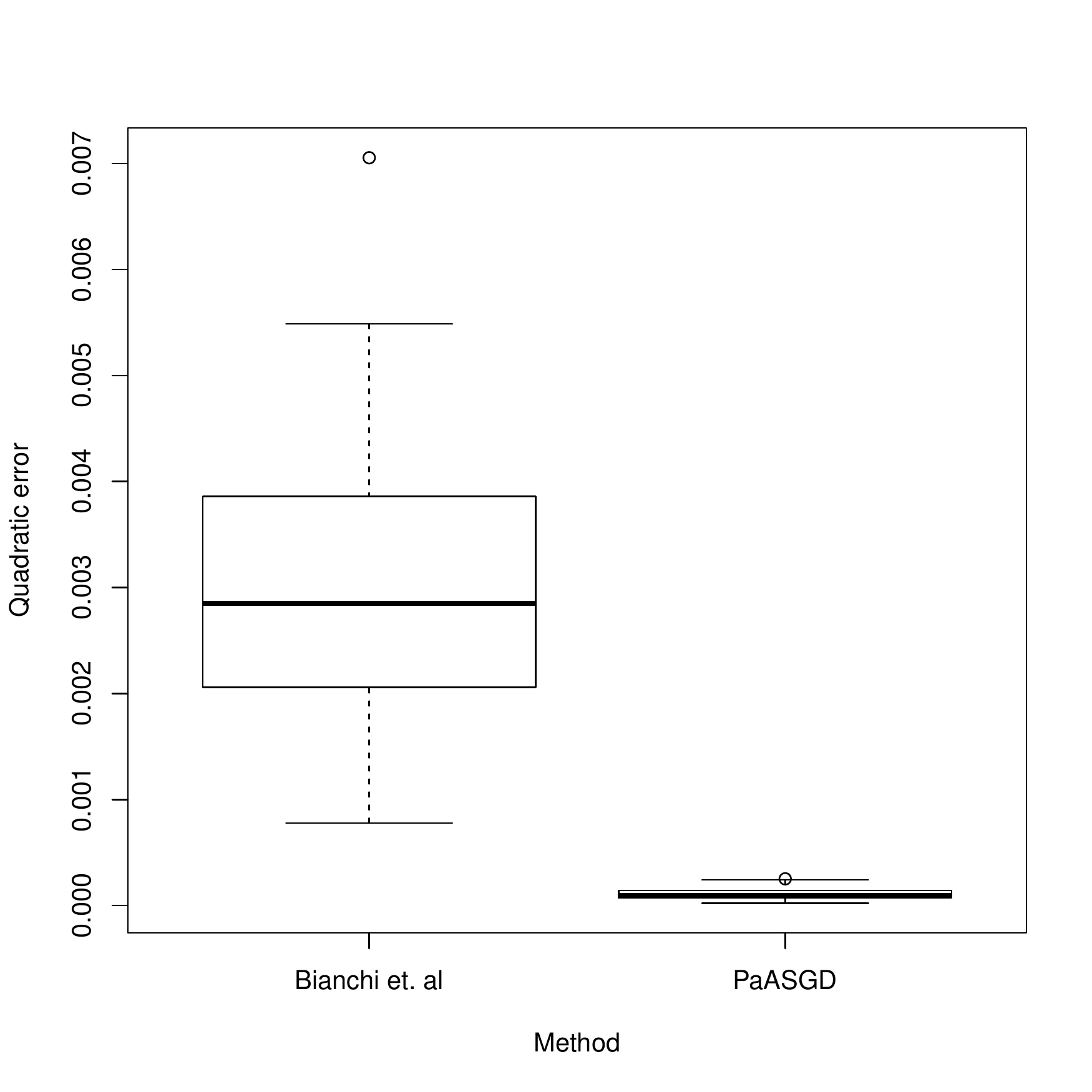}
\caption{Comparison between the quadratic errors obtained with the PASG algorithm and the ones obtained with the algorithm introduced by Bianchi et. al, for $p=10$ machines and with a Non-Equal Distribution.\label{figcomp}}
\end{center}
\end{figure}

In a recent work, \citep{bianchi2013performance} proposes to parallelize ASGD but for an uniform distribution of the data between the machines, and so that without taking the number of data per machine into  account. In Figure \ref{figcomp}, we show the significant improvement represented by our algorithm compare to the previous one. Indeed, although the algorithm proposed by \citep{bianchi2013performance} stay quite efficient in the context of Non-Equal Distribution, it can be very less accurate than the PSGD algorithm.

\bibliographystyle{plain}

\section{Proofs}\label{secproof}
\subsection{Some decompositions of the algorithms}
We first recall some usual decompositions of the algorithms, which will be useful in the proofs. First, for all $k \geq 1$, let us introduce the sequences $\left( \xi_{i,k} \right)$, defined, for all $i = 1,...,p$ and for all $k \geq 1$, by $\xi_{i,k+1} := \nabla G \left( m_{i,k} \right) - \nabla_{h} g \left( X_{i,k+1},m_{i,k}\right)$. Moreover, let us introduce the sequences of $\sigma$-algebras $\left( \mathcal{F}_{1,k} \right),..., \left( \mathcal{F}_{p,k} \right)$ defined for all $i=1,...,p$ and $k \geq 1$ by $\mathcal{F}_{i,k} := \sigma \left( X_{i,1},...,X_{i,k} \right)$. Then, for all $i=1,...,p$, $\left( \xi_{i,k} \right)_{k}$ is a sequence of martingale differences adapted to the filtration $\left( \mathcal{F}_{i,k} \right)_{k}$. Moreover, the SGD can be written, for all $i=1,...,p$ and for all $k \geq 1$, as
\begin{equation}
\label{decxi} m_{i,k+1} = m_{i,k} - \gamma_{k} \nabla G \left( m_{i,k} \right) + \gamma_{k}\xi_{i,k+1}.
\end{equation}
Moreover, linearizing the gradient, the SGD can be decomposed as
\begin{equation}\label{decdelta}
m_{i,k+1} - m = \left( I_{H} -\gamma_{k}\Gamma_{m}\right) \left( m_{i,k} - m \right) + \gamma_{k}\xi_{i,k+1} - \gamma_{k}\delta_{i,k},
\end{equation}
where $\delta_{i,k} := \nabla G \left( m_{i,k} \right) - \Gamma_{m} \left( m_{i,k} - m \right)$ is the remainder term in the Taylor's expansion of the gradient. Finally, summing these equalities, applying an Abel's transform and dividing by $n_{i}$, one can obtain (see \citep{Pel00})
\begin{align}
\Gamma_m\left[\notag \overline{m}_{i,n_{i}} - m\right] & = \frac{m_{i,1}-m}{n_{i}\gamma_{1}} - \frac{m_{i,n_{i}+1}-m}{n_{i}\gamma_{n_{i}}} + \frac{1}{n_{i}}\sum_{k=2}^{n_{i}}\left( m_{i,k} - m \right) \left( \frac{1}{\gamma_{k}} - \frac{1}{\gamma_{k-1}} \right) \\
&  - \frac{1}{n_{i}}\sum_{k=1}^{n_{i}}\delta_{i,k} + \frac{1}{n_{i}}\sum_{k=1}^{n} \xi_{i,k+1} .
\end{align}
Finally, by linearity, the PASG algorithm can be written as
\begin{align}\label{decpar}
\notag \Gamma_{m} \left[ \widehat{m}_{n} - m \right] & = \frac{1}{\sum_{i=1}^{p}n_{i}}\sum_{i=1}^{p}\frac{m_{i,1}-m}{\gamma_{1}} - \frac{1}{\sum_{i=1}^{p}n_{i}}\sum_{i=1}^{p}\frac{m_{i,n_{i}+1}-m}{\gamma_{n_{i}}} + \frac{1}{\sum_{i=1}^{p}n_{i}}\sum_{i=1}^{p}\sum_{k=2}^{n_{i}}\left( m_{i,k} - m \right) \left( \frac{1}{\gamma_{k}} - \frac{1}{\gamma_{k-1}} \right) \\
&  - \frac{1}{\sum_{i=1}^{p}n_{i}}\sum_{i=1}^{p}\sum_{k=1}^{n_{i}}\delta_{i,k} + \frac{1}{\sum_{i=1}^{p}n_{i}}\sum_{i=1}^{p}\sum_{k=1}^{n} \xi_{i,k+1} .
\end{align}

\subsection{Proof of Theorems \ref{theostr} and \ref{theoloc}}
In order to prove Theorems \ref{theostr} and \ref{theoloc}, we just have to bound each term on the right-hand side of (\ref{decpar}). 

\medskip

\noindent\textbf{Bounding $\mathbb{E}\left[ \left\| \frac{1}{\sum_{i=1}^{p}n_{i}}\sum_{i=1}^{p}\frac{m_{i,1}-m}{\gamma_{1}} \right\|^{2} \right]$.} With the help of Cauchy-Schwarz inequality and of inequality (\ref{vitlp}) or (\ref{vitl2}), one can check that 
\begin{equation}\label{maj1}
\mathbb{E}\left[ \left\| \frac{1}{\sum_{i=1}^{p}n_{i}}\sum_{i=1}^{p}\frac{m_{i,1}-m}{\gamma_{1}} \right\|^{2} \right] \leq \frac{p}{\left( \sum_{i=1}^{p}n_{i} \right)^{2}} \sum_{i=1}^{p} \mathbb{E}\left[ \left\| \frac{m_{i,1}-m}{\gamma_{1}}\right\|^{2}\right] \leq \frac{p^{2}C_{1}c_{\gamma}^{-2}}{\left(\sum_{i=1}^{p}n_{i}\right)^{2}} := A_{1,p}^2.
\end{equation}

\medskip

\noindent\textbf{Bounding $ \mathbb{E}\left[ \left\| \frac{1}{\sum_{i=1}^{p}n_{i}}\sum_{i=1}^{p}\frac{m_{i,n_{i}+1}-m}{\gamma_{n_{i}}} \right\|^{2} \right]$.} In the same way, with the help of Cauchy-Schwarz inequality and of inequality (\ref{vitlp}) or (\ref{vitl2}),
\begin{align*}
\mathbb{E}\left[ \left\| \frac{1}{\sum_{i=1}^{p}n_{i}}\sum_{i=1}^{p}\frac{m_{i,n_{i}+1}-m}{\gamma_{n_{i}}} \right\|^{2} \right] & \leq \frac{p}{\left( \sum_{i=1}^{p}n_{i}\right)^{2}}\sum_{i=1}^{p}\mathbb{E}\left[ \left\| \frac{m_{i,n_{i}+1}-m}{\gamma_{n_{i}}}  \right\|^{2} \right] \\
& \leq \frac{p}{\left( \sum_{i=1}^{p}n_{i}\right)^{2}}\sum_{i=1}^{p} C_{1}c_{\gamma}^{-2} (\frac{n_{i}}{n_i +1})^{\alpha}.
\end{align*}
This yields,
\begin{equation}\label{maj2}
\mathbb{E}\left[ \left\| \frac{1}{\sum_{i=1}^{p}n_{i}}\sum_{i=1}^{p}\frac{m_{i,n_{i}+1}-m}{\gamma_{n_{i}}} \right\|^{2} \right]  \leq \frac{p^{2}c_{\gamma}^{-2}C_{1}}{\left( \sum_{i=1}^{p}n_{i}\right)^{2}}:= A_{2,p}^2.
\end{equation}

\medskip

\noindent \textbf{Bounding $\frac{1}{\sum_{i=1}^{p}n_{i}}\sum_{i=1}^{p}\sum_{k=2}^{n_{i}}\left( m_{i,k} - m \right) \left( \frac{1}{\gamma_{k}} - \frac{1}{\gamma_{k-1}} \right)$.} In the same way, applying Lemma 4.3 in~\citep{godichon2015},
\begin{align*}
\mathbb{E}\left[ \left\| \frac{1}{\sum_{i=1}^{p}n_{i}}\sum_{i=1}^{p}\sum_{k=2}^{n_{i}}\left( m_{i,k} - m \right) \left( \frac{1}{\gamma_{k}} - \frac{1}{\gamma_{k-1}} \right) \right\|^{2} \right] & \leq \frac{p}{\left( \sum_{i=1}^{p}n_{i}\right)^{2}}\sum_{i=1}^{p}\mathbb{E}\left[ \left\| \sum_{k=2}^{n_{i}}\left( m_{i,k} - m \right) \left( \frac{1}{\gamma_{k}} - \frac{1}{\gamma_{k-1}} \right) \right\|^{2} \right] \\
& \leq \frac{p}{\left( \sum_{i=1}^{p}n_{i}\right)^{2}}\sum_{i=1}^{p}\left( \sum_{k=2}^{n_{i}}\sqrt{\mathbb{E}\left[ \left\| m_{i,k} - m \right\|^{2} \right]} \left| \frac{1}{\gamma_{k}} - \frac{1}{\gamma_{k-1}} \right| \right)^{2}. 
\end{align*}
Since 
\[\left| \frac{1}{\gamma_{k}}- \frac{1}{\gamma_{k-1}}\right| ~\leq ~\alpha ~c_{\gamma}^{-1}(k-1)^{\alpha-1}, 
\]
and applying inequality (\ref{vitlp}) or (\ref{vitl2})
\begin{align*}
\mathbb{E}\left[ \left\| \frac{1}{\sum_{i=1}^{p}n_{i}}\sum_{i=1}^{p}\sum_{k=2}^{n_{i}}\left( m_{i,k} - m \right) \left( \frac{1}{\gamma_{k}} - \frac{1}{\gamma_{k-1}} \right) \right\|^{2} \right] 
& \leq \frac{p\alpha^{2} c_{\gamma}^{-2}C_{1}}{\left( \sum_{i=1}^{p}n_{i}\right)^{2}}\sum_{i=1}^{p}\left( \sum_{k=2}^{n_{i}}\frac{(k-1)^{\alpha-1}}{k^{\alpha/2}} \right)^{2}\\
&\leq  \frac{p\alpha^{2} c_{\gamma}^{-2}C_{1}}{\left( \sum_{i=1}^{p}n_{i}\right)^{2}}\sum_{i=1}^{p}\left( \sum_{k=2}^{n_{i}}\frac{(k-1)^{\alpha-1}}{(k-1)^{\alpha/2}} \right)^{2}\\
&\leq \frac{p\alpha^{2} c_{\gamma}^{-2}C_{1}}{\left( \sum_{i=1}^{p}n_{i}\right)^{2}}\sum_{i=1}^{p}\left( \sum_{k=2}^{n_{i}}(k-1)^{\alpha/2-1} \right)^{2}.
\end{align*}
Then, since $\alpha <1$, with the help of an integral test for convergence and thanks to Hölder's inequality,
\begin{align}
\notag \mathbb{E}\left[ \left\| \frac{1}{\sum_{i=1}^{p}n_{i}}\sum_{i=1}^{p}\sum_{k=2}^{n_{i}}\left( m_{i,k} - m \right) \left( \frac{1}{\gamma_{k}} - \frac{1}{\gamma_{k-1}} \right) \right\|^{2} \right] & \leq \frac{4p\alpha c_{\gamma}^{-2}C_{1}}{\left( \sum_{i=1}^{p}n_{i}\right)^{2}}\sum_{i=1}^{p} (n_{i}-1)^{\alpha} \\
\notag & \leq \frac{4p^{2-\alpha} \alpha c_{\gamma}^{-2}C_{1}}{\left( \sum_{i=1}^{p}n_{i}\right)^{2}}\left( \sum_{i=1}^{p}(n_{i}-1) \right)^{\alpha} \\
\label{maj3}& \leq \frac{4p^{2-\alpha} \alpha c_{\gamma}^{-2}C_{1}}{\left( \sum_{i=1}^{p}n_{i}\right)^{2-\alpha}}:= A_{3,p}^2.
\end{align}
\medskip

\noindent\textbf{Bounding $\mathbb{E}\left[ \left\| \frac{1}{\sum_{i=1}^{p}n_{i}}\sum_{i=1}^{p}\sum_{k=1}^{n_{i}}\delta_{i,k} \right\|^{2}\right]$.} First, let us recall that  there is a positive constant $C_{m}$ such that for all $i=1,...,p$, and for all integer $k$,
\begin{equation}
\label{majdelta} \left\| \delta_{i,k} \right\| \leq C_{m} \left\| m_{i,k} - m \right\|^{2}.
\end{equation}
Moreover, thanks to Lemma 4.1 in \citep{godichon2015},
\begin{align*}
\mathbb{E}\left[ \left\| \frac{1}{\sum_{i=1}^{p}n_{i}}\sum_{i=1}^{p}\sum_{k=1}^{n_{i}}\delta_{i,k} \right\|^{2}\right] & \leq \frac{p}{\left( \sum_{i=1}^{p}n_{i}\right)^{2}}\sum_{i=1}^{p} \mathbb{E}\left[ \left\| \sum_{k=1}^{n_{i}} \delta_{i,k} \right\|^{2} \right] \leq \frac{p}{\left( \sum_{i=1}^{p}n_{i}\right)^{2}}\sum_{i=1}^{p} \left( \sum_{k=1}^{n_{i}} \sqrt{\mathbb{E}\left[ \left\| \delta_{i,k}\right\|^{2}\right]} \right)^{2}.
\end{align*}
Then, applying inequalities (\ref{majdelta}) and (\ref{vitlp}) or (\ref{vitl4}),
\begin{align*}
\mathbb{E}\left[ \left\| \frac{1}{\sum_{i=1}^{p}n_{i}}\sum_{i=1}^{p}\sum_{k=1}^{n_{i}}\delta_{i,k} \right\|^{2}\right] & \leq \frac{C_{m}^{2}p}{\left( \sum_{i=1}^{p}n_{i}\right)^{2}}\sum_{i=1}^{p} \left( \sum_{k=1}^{n_{i}} \sqrt{\mathbb{E}\left[ \left\| m_{i,k} - m\right\|^{4}\right]} \right)^{2} 
\leq \frac{C_{m}^{2}C_{2}p}{\left( \sum_{i=1}^{p}n_{i}\right)^{2}}\sum_{i=1}^{p} \left( \sum_{k=1}^{n_{i}} k^{-\alpha} \right)^{2}.
\end{align*}
Thus, since $1/2<\alpha<1$, with the help of an integral test for convergence and thanks to Hölder's inequality,
\begin{align}
\notag \mathbb{E}\left[ \left\| \frac{1}{\sum_{i=1}^{p}n_{i}}\sum_{i=1}^{p}\sum_{k=1}^{n_{i}}\delta_{i,k} \right\|^{2}\right] & \leq \frac{C_{m}^{2}C_{2}(1-\alpha)^{-2}p}{\left( \sum_{i=1}^{p}n_{i}\right)^{2}}\sum_{i=1}^{p} n_{i}^{2-2\alpha} \leq \frac{C_{m}^{2}C_{2}(1-\alpha)^{-2}p^{2\alpha}}{\left( \sum_{i=1}^{p}n_{i}\right)^{2}}\left( \sum_{i=1}^{p}n_{i} \right)^{2-2\alpha} \\
\label{maj4}&\leq \frac{C_{m}^{2}C_{2}(1-\alpha)^{-2}p^{2\alpha}}{(\sum_{i=1}^{p}n_{i})^{2\alpha}}:= A_{4,p}^2.
\end{align}
\medskip

\noindent\textbf{Bounding $\mathbb{E}\left[ \left\| \frac{1}{\sum_{i=1}^{p}n_{i}}\sum_{i=1}^{p}\sum_{k=1}^{n} \xi_{i,k+1}\right\|^{2} \right]$.} First, by definition of the sequences $\left( \xi_{i,k} \right)$  and thanks to assumption \textbf{(H4)} or \textbf{(H8)}, for all $i=1,...,p$ and for all positive integer $k$,
\begin{align}
\notag  \mathbb{E}\left[ \left\| \xi_{i,k+1} \right\|^{2} \right] & = \mathbb{E}\left[ \left\| \nabla_{h}g \left( X_{i,+1},m_{i,k} \right) \right\|^{2}\right] - \mathbb{E}\left[ \left\| \nabla G \left( m_{i,k} \right) \right\|^{2} \right] \\
& \leq L_{1} + L_{1} \mathbb{E}\left[ \left\| m_{i,k} - m \right\|^{2} \right] . \label{majxi}
\end{align}
Moreover, since for all $i=1,...,p$, $\left( \xi_{i,k} \right)$ is a sequence of martingale differences adapted to the filtration $\left( \mathcal{F}_{i,k} \right)$ and since for all $i=1,...,p$ and $j=1,...,p$ such taht $i \neq j$ the sequences $\left( \xi_{i,k}\right)$ and $\left( \xi_{j,k}\right)$ are independent,
\begin{align*}
& \mathbb{E}\left[ \sum_{k=1}^{n_{i}}\xi_{i,k} \right] = 0 \quad \forall i=1,...,p, \quad \text{and} \quad \mathbb{E}\left[ \left\langle \sum_{k=1}^{n_{i}}\xi_{i,k} , \sum_{k=1}^{n_{j}}\xi_{j,k} \right\rangle \right] = 0, \quad \forall i,j =1,...,p \quad \text{s.t} \quad i\neq j.
\end{align*}
Then, applying inequality (\ref{majxi}),
\begin{align*}
\mathbb{E}\left[ \left\| \frac{1}{\sum_{i=1}^{p}n_{i}}\sum_{i=1}^{p}\sum_{k=1}^{n_{i}} \xi_{i,k+1}\right\|^{2} \right] & = \frac{1}{\left( \sum_{i=1}^{p}n_{i} \right)^{2}}\sum_{i=1}^{p} \mathbb{E}\left[ \left\| \sum_{k=1}^{n_{i}} \xi_{i,k+1}\right\|^{2}\right] \\
& = \frac{1}{\left( \sum_{i=1}^{p}n_{i} \right)^{2}}\sum_{i=1}^{p}  \sum_{k=1}^{n_{i}}\mathbb{E}\left[ \left\| \xi_{i,k+1}\right\|^{2}\right] \\
& \leq \frac{1}{\left( \sum_{i=1}^{p}n_{i} \right)^{2}}\sum_{i=1}^{p}  \sum_{k=1}^{n_{i}}\left( L_{1} + L_{1}\mathbb{E}\left[ \left\| m_{i,k} - m \right\|^{2} \right]\right)
\end{align*}
Then, thanks to inequality (\ref{vitlp}) or (\ref{vitl2}), and with the help of an integral test for convergence,
\begin{align*}
\mathbb{E}\left[ \left\| \frac{1}{\sum_{i=1}^{p}n_{i}}\sum_{i=1}^{p}\sum_{k=1}^{n_{i}} \xi_{i,k+1}\right\|^{2} \right] & \leq \frac{L_{1}}{\sum_{i=1}^{p}n_{i}}+  \frac{L_{1}C_{1}}{\left( \sum_{i=1}^{p}n_{i} \right)^{2}}\sum_{i=1}^{p}\sum_{k=1}^{n_{i}} k^{-\alpha} \\
& \leq \frac{L_{1}}{\sum_{i=1}^{p}n_{i}}+  \frac{L_{1}C_{1}(1-\alpha )^{-1}}{\left( \sum_{i=1}^{p}n_{i} \right)^{2}}\sum_{i=1}^{p}n_{i}^{1-\alpha}.
\end{align*}
Finally, applying Hölder's inequality and since $\alpha < 1$,
\begin{align}
\label{maj5}\mathbb{E}\left[ \left\| \frac{1}{\sum_{i=1}^{p}n_{i}}\sum_{i=1}^{p}\sum_{k=1}^{n_{i}} \xi_{i,k+1}\right\|^{2} \right] & \leq \frac{L_{1}}{\sum_{i=1}^{p}n_{i}}+  \frac{L_{1}C_{1}(1-\alpha )^{-1}p^{\alpha}}{\left( \sum_{i=1}^{p}n_{i} \right)^{1+\alpha}}:= A_{6,p}^2 + A_{5,p}^2.
\end{align}
\medskip

\noindent\textbf{Conclusion.} Since the smallest eigenvalue (or the limit inf of the eignevalues for infinite dimensional spaces) of $\Gamma_{m}$ denoted by $\lambda_{\min}$ is positive, let $n = \sum_{i=1}^{p}n_{i}$,
\begin{align*}
\mathbb{E}\left[ \left\| \widehat{m}_{n} - m \right\|^{2} \right] \leq \frac{1}{\lambda_{\min}^{2}} \mathbb{E}\left[ \left\| \Gamma_{m}\left( \widehat{m}_{n} - m \right) \right\|^{2} \right].
\end{align*}
Then, applying Cauchy-Schwarz's inequality as well as inequalities (\ref{maj1}) to (\ref{maj5}), one can check that
\[
\mathbb{E}\left[ \left\| \widehat{m}_{n} - m \right\|^{2} \right] \leq \sum_{j=1}^{6}\lambda_{\min}^{-2}A_{j,p,n}^2 + \sum_{j,j'=1,\quad j\neq j'}^{6}\lambda_{\min}^{-2}A_{j,p,n}A_{j',p,n} ,
\]
which concludes the proof.

\subsection{Proof of Theorem \ref{tlcstr} and \ref{tlcloc}}
First, one can check that the first term on the right-hand side of equality (\ref{decpar}) are negligeable, i.e
\begin{align*}
& \frac{1}{\sqrt{\sum_{i=1}^{p}n_{i}}}\sum_{i=1}^{p}\frac{m_{i,1}-m}{\gamma_{1}} = o \left( 1 \right) \quad \quad \mathbb{P}, \\
&  \frac{1}{\sqrt{\sum_{i=1}^{p}n_{i}}}\sum_{i=1}^{p}\frac{m_{i,n_{i}+1}-m}{\gamma_{n_{i}}} = o \left( 1 \right) \quad \quad \mathbb{P},  \\
& \frac{1}{\sqrt{\sum_{i=1}^{p}n_{i}}}\sum_{i=1}^{p}\sum_{k=2}^{n_{i}}\left( m_{i,k} - m \right) \left( \frac{1}{\gamma_{k+1}} - \frac{1}{\gamma_{k}} \right)  = o \left( 1 \right) \quad \quad \mathbb{P}, \\
&   \frac{1}{\sqrt{\sum_{i=1}^{p}n_{i}}}\sum_{i=1}^{p}\sum_{k=1}^{n_{i}}\delta_{i,k} = o \left( 1 \right) \quad \quad \mathbb{P}.
\end{align*}
Indeed, it is a direct application of inequalities (\ref{maj1}) to (\ref{maj5}) when assumptions \textbf{(H1)} to \textbf{(H5)} are verified and a direct application of Theorem 4.1 in \citep{GB2016} when assumptions \textbf{(H1)} to \textbf{(H4)} and \textbf{(H7)} are verified. In order to get the asymptotic normality of the term $\left( \sum_{i=1}^{p}\sum_{k=1}^{n_{i}}\xi_{i,k+1} \right)$, let us first remark that with a good choice of index, this term can be seen as a sum of martingale differences term. Then, we just have to check that assumptions of Theorem 5.1 in \citep{Jak88} are fulfilled, i.e let $\left( e_{j} \right)_{j \in J}$ be an orthonormal basis of $H$ and $\psi_{j,j'} := \left\langle \Sigma e_{j},e_{j'} \right\rangle$ for all $j,j' \in J$, we have to verify:
\begin{equation}\label{eq1}
\forall \eta > 0, \quad \lim_{n \to \infty} \mathbb{P} \left( \sup_{1\leq k \leq n_{i},i=1,\dots p} \frac{1}{\sqrt{\sum_{i=1}^{p}n_{i}}}\left\| \xi_{i,k+1} \right\| > \eta \right) =0,
\end{equation}
\begin{equation}\label{eq2}
\lim_{n \to \infty} \frac{1}{\sum_{i=1}^{p}n_{i}}\sum_{i=1}^{p}\sum_{k=1}^{n_{i}} \left\langle \xi_{i,k+1},e_{j}\right\rangle \left\langle \xi_{i,k+1},e_{j'}\right\rangle = \psi_{j,j'} \quad a.s, \quad \forall j,j' \in J,
\end{equation}
\begin{equation}\label{eq3}
\forall \epsilon > 0 , \quad \lim_{N \to \infty} \limsup_{n \to \infty} \mathbb{P}\left( \frac{1}{\sum_{i=1}^{p}\sum_{k=1}^{n_{i}}}\sum_{j=N}^{\infty} \left\langle \xi_{i,k+1},e_{j} \right\rangle^{2}> \epsilon \right) = 0. 
\end{equation}

\medskip

\noindent\textbf{Proof of (\ref{eq1}):} Let $\eta > 0$, applying Markov's inequality,
\begin{align*}
\mathbb{P} \left( \sup_{1\leq k \leq n_{i},i=1,\dots p} \frac{1}{\sqrt{\sum_{i=1}^{p}n_{i}}}\left\| \xi_{i,k+1} \right\| > \eta \right) & = \sum_{i=1}^{p}\sum_{k=1}^{n_{i}} \mathbb{P}\left( \frac{1}{\sqrt{\sum_{i=1}^{p}n_{i}}}\left\| \xi_{k+1} \right\| > \eta \right) \\
& \leq \frac{\eta^{-4}}{\left( \sum_{i=1}^{p}n_{i} \right)^{2}}\sum_{i=1}^{p}\sum_{k=1}^{n_{i}} \mathbb{E}\left[ \left\| \xi_{i,k+1} \right\|^{4} \right] .
\end{align*}
Moreover, note that thanks to Assumption \textbf{(H4)} or \textbf{(H8)}, for all $i=1,...,p$ and $1 \leq k \leq n_{i}$,
\begin{align*}
\mathbb{E}\left[ \left\| \xi_{i,k+1} \right\|^{4} \right] \leq 2^{4} \mathbb{E}\left[ \left\| \nabla_{h} g \left( X_{i,k+1},m_{i,k} \right) \right\|^{4} \right] \leq 2^{4}L_{2} \left( 1+ \mathbb{E}\left[ \left\| m_{i,k} \right\|^{4} \right] \right) \leq 2^{4}L_{2} \left( 1 + C_{2} \right) .
\end{align*}
Then, 
\begin{align*}
\mathbb{P} \left( \sup_{1\leq k \leq n_{i},i=1,\dots p} \frac{1}{\sqrt{\sum_{i=1}^{p}n_{i}}}\left\| \xi_{i,k+1} \right\| > \eta \right) \leq \frac{\eta^{-4}2^{4}L_{2}\left( 1+C_{2}\right)}{\sum_{i=1}^{p}n_{i}}.
\end{align*}

\medskip

\textbf{Proof of (\ref{eq2}):} First, let $\otimes$ be the bilinear application defined for all $h,h',h'' \in H$ by $\left( h\otimes h' \right)(h'') = \left\langle h,h'' \right\rangle h'$. Note that
\[
\frac{1}{\sum_{i=1}^{p}n_{i}}\sum_{i=1}^{p}\sum_{k=1}^{n_{i}}\xi_{i,+1} \otimes \xi_{i,k+1} = \frac{1}{\sum_{i=1}^{p}n_{i}}\sum_{i=1}^{p}\sum_{k=1}^{n_{i}} \mathbb{E}\left[ \xi_{i,k+1} \otimes \xi_{i,k+1} |\mathcal{F}_{i,k} \right] + \frac{1}{\sum_{i=1}^{p}n_{i}}\sum_{i=1}^{p}\sum_{k=1}^{n_{i}} \epsilon_{i,k+1}, 
\]
with $\epsilon_{i,k+1} := \xi_{i,k+1} \otimes \xi_{i,k+1} - \mathbb{E}\left[ \xi_{i,k+1} \otimes \xi_{i,k+1} |\mathcal{F}_{i,k} \right]$. Remark that we a good choice of index, $\left( \epsilon_{i,k} \right)_{i,k}$ can be seen as a sequence of martingale differences, and one can check that
\[
\lim_{n \to \infty} \frac{1}{\sum_{i=1}^{p}n_{i}} \sum_{i=1}^{p}\sum_{k=1}^{n_{i}}\epsilon_{i,k+1} = 0 \quad a.s.
\] 
Let us now prove that the sequence of operators $\left( \mathbb{E}\left[ \xi_{i,k+1} \otimes \xi_{i,k+1} |\mathcal{F}_{i,k} \right] \right)_{i,k}$ converges almost surely to $\Sigma$ with respect to the Frobenius norm when $k$ goes to infinity. First, note that
\begin{align*}
& \left\| \mathbb{E}\left[ \xi_{i,k+1} \otimes \xi_{i,k+1} |\mathcal{F}_{k} \right] - \Sigma \right\|_{F} \\
&  = \left\| \mathbb{E}\left[ \nabla_{h}g \left( X_{i,k+1},m_{i,k}\right) \otimes \nabla_{h}g \left( X_{i,k+1},m_{i,k}\right) |\mathcal{F}_{i,k} \right] - \Sigma - \nabla G \left( m_{i,k} \right) \otimes \nabla G \left( m_{i,k} \right)  \right\|_{F} \\
& \leq \left\| \mathbb{E}\left[ \nabla_{h}g \left( X_{i,k+1},m_{i,k}\right) \otimes \nabla_{h}g \left( X_{i,k+1},m_{i,k}\right) |\mathcal{F}_{i,k} \right] - \Sigma  \right\|_{F} + \left\|  \nabla G \left( m_{i,k} \right) \otimes \nabla G \left( m_{i,k} \right)  \right\|_{F} .
\end{align*}
Then, thanks to Assumption \textbf{(H6)}, since $\left\| \nabla G \left( m_{i,k} \right) \right\| \leq C \left\| m_{i,k} - m \right\|$ for all $i,k$ (see \citep{GB2016}), and since $\left( m_{i,k} \right)$ converges almost surely to $m$,
\begin{align*}
& \lim_{k \to \infty} \left\| \mathbb{E}\left[ \nabla_{h}g \left( X_{i,k+1},m_{i,k}\right) \otimes \nabla_{h}g \left( X_{i,k+1},m_{i,k}\right) |\mathcal{F}_{i,k} \right] - \Sigma  \right\|_{F} = 0 \quad a.s, \quad \forall i=1,...,p, \\
& \lim_{k \to \infty}\left\|  \nabla G \left( m_{i,k} \right) \otimes \nabla G \left( m_{i,k} \right)  \right\|_{F} = \lim_{k\to \infty} \left\| \nabla G \left( m_{i,k} \right) \right\|_{F}^{2} = 0 \quad a.s, \quad \forall i=1,...,p.
\end{align*}
Then, the sequences $\left( \mathbb{E}\left[ \xi_{i,k+1} \otimes \xi_{i,k+1} |\mathcal{F}_{k} \right] \right)_{k \geq 1}$ converges almost surely to $\Sigma$ with respect to the Frobenius norm and as a consequence, for all $j,j' \in J$, 
\[
\lim_{k \to \infty} \left\langle \mathbb{E}\left( \xi_{i,k+1} \otimes \xi_{i,k+1} |\mathcal{F}_{i,k} \right] (e_{j}) , e_{j'} \right\rangle = \psi_{j,j'} \quad a.s, \quad \forall i=1,...,p.
\]
Thus, applying Toeplitz's lemma, for all $j,j' \in J$,
\[
\lim_{n\to \infty} \frac{1}{\sum_{i=1}^{p}n_{i}}\sum_{i=1}^{p}\sum_{k=1}^{n_{i}} \left\langle \mathbb{E}\left( \xi_{i,k+1} \otimes \xi_{i,k+1} |\mathcal{F}_{i,k} \right] (e_{j}) , e_{j'} \right\rangle = \psi_{j,j'} \quad a.s.
\]

\medskip

\noindent\textbf{Proof of (\ref{eq3}):} Let $\epsilon > 0$, applying Markov's inequality,
\begin{align*}
\mathbb{P} \left( \frac{1}{\sum_{i=1}^{p}n_{i}} \sum_{i=1}^{p}\sum_{k=1}^{n_{i}}\sum_{j=N}^{\infty} \left\langle \xi_{i,k+1} , e_{j} \right\rangle > \epsilon \right) & \leq \frac{\epsilon^{-2}}{\sum_{i=1}^{p}n_{i}} \sum_{i=1}^{p}\sum_{k=1}^{n_{i}} \sum_{j=N}^{\infty} \mathbb{E}\left[ \left\langle \xi_{i,k+1},e_{j} \right\rangle^{2} \right] \\
& \leq  \frac{\epsilon^{-2}}{\sum_{i=1}^{p}n_{i}} \sum_{i=1}^{p}\sum_{k=1}^{n_{i}}\sum_{j=N}^{\infty} \mathbb{E}\left[ \mathbb{E}\left[  \left\langle \xi_{i,k+1},e_{j} \right\rangle^{2} |\mathcal{F}_{i,k} \right] \right] .
\end{align*}
Since for all $j \in J$, $\left\langle \xi_{i,k+1},e_{j} \right\rangle^{2} = \left\langle \xi_{i,k+1} \otimes \xi_{i,k+1} \left( e_{j} \right) , e_{j} \right\rangle$, by linearity and by dominated convergence,
\begin{align*}
\mathbb{P} \left( \frac{1}{\sum_{i=1}^{p}n_{i}} \sum_{i=1}^{p}\sum_{k=1}^{n_{i}}\sum_{j=N}^{\infty} \left\langle \xi_{i,k+1} , e_{j} \right\rangle > \epsilon \right) & \leq \frac{1}{\epsilon^{2}}\sum_{j=N}^{\infty} \frac{1}{\sum_{i=1}^{p}n_{i}}\sum_{i=1}^{p}\sum_{k=1}^{n_{i}} \mathbb{E}\left[ \mathbb{E}\left[  \left\langle \xi_{i,k+1} \otimes \xi_{i,k+1} \left( e_{j} \right) , e_{j} \right\rangle |\mathcal{F}_{i,k} \right] \right] \\
& = \frac{1}{\epsilon^{2}}\sum_{j=N}^{\infty} \frac{1}{\sum_{i=1}^{p}n_{i}}\sum_{i=1}^{p}\sum_{k=1}^{n_{i}} \mathbb{E}\left[   \left\langle \mathbb{E}\left[ \xi_{i,k+1} \otimes \xi_{i,k+1} |\mathcal{F}_{i,k} \right] \left( e_{j}\right)  , e_{j} \right\rangle  \right] .
\end{align*}
Since $\mathbb{E}\left[ \xi_{i,k+1} \otimes \xi_{i,k+1} |\mathcal{F}_{i,k} \right]$ converges almost surely to $\Sigma$ and by dominated convergence,
\[
\limsup_{n} \mathbb{P}\left( \frac{1}{n} \sum_{i=1}^{p}\sum_{k=1}^{n_{i}}\sum_{j=N}^{\infty} \left\langle \xi_{i,k+1},e_{j} \right\rangle > \epsilon \right) \leq \frac{1}{\epsilon^{2}} \sum_{j=N}^{\infty} \left\langle \Sigma (e_{j}),e_{j} \right\rangle,
\]
and one can conclude as in \citep{GB2017}.

\end{document}